\documentclass[11pt]{article}
\usepackage{amsmath}
\usepackage{amssymb}
\usepackage{amscd}
\usepackage[all]{xy}
\usepackage{epsfig}
\def\B{\mathcal B}
\def\qsl{\mathcal O(SL_q(2))}
\def\N{\mathbb N}
\def\C{\mathbb C}
\def\R{\mathbb R}

\newtheorem{theo}{Theorem}[section]
\newtheorem{prop}[theo]{Proposition}

\newtheorem{lemm}[theo]{Lemma}
\newtheorem{defi}[theo]{Definition}
\newtheorem{rem}[theo]{Remark}

\title{\textbf{The Representation Category of the Quantum Group
of a Non-degenerate Bilinear Form}} 
\author{Julien Bichon}
\date{{\small \textsl{Laboratoire de Math\'ematiques Appliqu\'ees,
Universit\'e de Pau et des Pays de l'Adour, \\
IPRA, Avenue de l'universit\'e,
64000 Pau, France.}
E-mail: Julien.Bichon@univ-pau.fr}}

\makeatletter
\renewcommand{\@makefnmark}{}
\makeatother

\begin{document}

\maketitle

\begin{abstract}
We show that the representation category of the quantum group of a 
non-degenerate bilinear form is monoidally equivalent
to the representation category of the quantum group $SL_q(2)$ for a 
well-chosen non-zero parameter $q$.
The key ingredient for the proof of this result
is the direct and explicit construction of an 
appropriate Hopf bigalois extension. 
Then we get, when the base field is of characteristic zero, a
full description of cosemisimple Hopf algebras whose representation semi-ring
is isomorphic to the one of $SL(2)$.
\end{abstract}

Keywords: Hopf algebra, monoidal category, Hopf-Galois extension.

\section{Introduction and main results}

Let $k$ be a commutative algebraically closed field, let $n \in \mathbb N^*$,
$n \geq 2$ and let $E \in GL(n)$. We consider the following algebra
$\mathcal B(E)$: it is the universal algebra with generators
$(a_{ij})_{1 \leq i,j \leq n}$ and satisfying the relations

$$E^{-\!1} {^t \! a} E a = I = a E^{-\!1} {^t \! a} E,$$
where $a$ is the matrix $(a_{ij})_{1 \leq i,j \leq n}$ and $I$ is the
identity matrix. This algebra admits a natural Hopf algebra
structure and was introduced by M. Dubois-Violette and G. Launer 
\cite{[DVL]}.
It is the function algebra on the
quantum (symmetry) group of a non-degenerate bilinear
form (see section 2). Let $q \in k^*$. For a well-chosen matrix 
$E_q \in GL(2)$, we have $\mathcal B(E_q) = \mathcal O(SL_q(2))$,
the function algebra on the quantum group $SL_q(2)$.
The main result of this paper describes the category of comodules
over $\mathcal B(E)$ for a general matrix $E$:

\begin{theo}
Let $E \in GL(n)$, $n \geq 2$, and let $q \in k^*$ be such that
$q^2 + {\rm tr}(E ^t \!  E^{-1})q +1 = 0$. 
Then we have an equivalence of monoidal categories:
$${\rm Comod}(\mathcal B(E)) \cong^{\otimes} 
{\rm Comod}(\mathcal O(SL_q(2))$$
between the comodule categories of $\mathcal B(E)$
and $\mathcal O( SL_q(2))$ respectively.
\end{theo}

When $k = \mathbb C$, T. Banica \cite{[Ba1]} proved a similar result,
in the compact quantum group case,
at the representation semi-ring level.
Theorem 1.1 covers the cosemisimple non-compact case as well
as the non-cosemisimple case. 
There are also other related results in the literature, in the
$SL(N)$ case: again by Banica \cite{[Ba2]} in the compact case 
and by Phung Ho Hai \cite{[Ph]}
in the cosemisimple case ($q$ is not a root of unity) in characteristic 
zero. In these two approaches the authors study Hopf algebras
reconstructed from  Hecke symmetries.

We wish to emphasize that our result is characteristic-free and does
not depend on the cosemisimplicity of the considered Hopf algebras.
The main reason is that our technique of proof is different from
the one of Banica and Phung Ho Hai. These two authors use reconstruction
techniques. Here we directly construct an explicit
$\qsl$-$\B(E)$-bigalois extension: by a very useful theorem of 
P. Schauenburg \cite{[Sc]} (see also K.H. Ulbrich \cite{[U]}), 
this is equivalent to construct an equivalence of monoidal 
categories between the comodule categories of these two Hopf algebras.
The technical difficulty in our approach is to show that that the
algebra we construct is non-zero. Since the the monoidal equivalences 
we get do not preserve the dimensions of the underlying
vector spaces in general, our  Galois extensions will be non-cleft in general.
The existence of non-cleft Hopf-Galois extensions was known: first by the
end of the paper \cite{[Br]} of A. Brugui\`eres, and also
by the results of Banica and Phung Ho Hai. However it is the first
time, at least to the best of our knowledge, that
non-cleft Hopf-Galois extensions are explicitly described.

Let us point out a negative consequence of Theorem 1.1 in the 
perspective of knot theory. Recall that the Jones polynomial
may be constructed from the representation category of the quantum
group $SL_q(2)$ (see the book \cite{[Ka]}). Theorem 1.1
means that one cannot expect to get any new link invariant from the
more general Hopf algebras $\B(E)$. 

\medskip

We also prove a kind of converse to Theorem 1.1: the description
of all cosemisimple Hopf algebras whose representation semi-ring is 
isomorphic to the one of $SL(2)$. 
Here we have to assume that the characteristic of
$k$ is zero. We say that an element $q \in k^*$ is generic if 
$q \in \{\pm 1 \}$ or if $q$ is not a root of unity.

\begin{theo} Let $A$ be a cosemisimple Hopf algebra whose representation 
semi-ring is isomorphic to the one of $SL(2)$. Then there exists 
$E \in GL(n)$ ($n\geq 2$) such that $A$ is isomorphic with 
$\mathcal B(E)$, and such that 
any solution of the equation
$q^2 + {\rm tr}(E ^t \!  E^{-1})q +1 = 0$ is generic.
If $F \in GL(m)$ is another matrix such that $A$ is isomorphic
with $\mathcal B(F)$, then $n=m$ and there exists $M \in GL(m)$
such that $F = {^t \! M}EM$
\end{theo} 

Once again an analogue of Theorem 1.2 was proved in the compact quantum group
case in \cite{[Ba1]}. But again we have here the cosemisimple non-compact case.
Theorem 1.2 was already known if one requires the fundamental comodule
of $A$ to be of dimension 2, partially by  results of S.L. Woronowicz
\cite{[Wo]}, a complete proof being given in P. Podle\'s and 
E. M\"{u}ller's 
notes \cite{[PM]}. The $SL(3)$-case has  been done by C. Ohn \cite{[O]}
with a constraint on the dimension of the fundamental comodule.
Finally the compact case $SU(N)$ was done in \cite{[Ba2]},
without any dimension constraint but
without an isomorphic classification.

Theorem 1.1 is used in an essential way to prove Theorem 1.2. 
The other main ingredient for the proof of Theorem 1.2 is the
representation theory of $SL_q(2)$, including the root of unity case
(see \cite{[KP]}). The  strategy of proof is then the same as the
one of Podle\'s and M\"{u}ller \cite{[PM]}. 

\medskip

Our paper is organized as follows. In Section 2 we briefly recall
some facts concerning the Hopf algebras $\B(E)$. In Section 3,
we associate an algebra $\B(E,F)$ to each pair $(E,F)$ of matrices.
It is shown that if $\B(E,F)$ is a non zero-algebra, then 
$\B(E,F)$ is a $\B(E)$-$\B(F)$-bigalois extension. 
In section 4 we prove, using the diamond Lemma \cite{[Be]},
 that $\B(E_q,F)$ is a non-zero algebra
for a well chosen $q \in k^*$: this proves Theorem 1.1. 
In section 5 we prove Theorem 1.2
and describe the isomorphic classification 
of the Hopf algebras $\B(E)$ (in characteristic zero). 
Finally we study possible
CQG algebra structures on $\B(E)$ in Section 6.

\smallskip

Throughout this paper $k$ is an algebraically closed field.

\section{The Hopf algebras $\mathcal B(E)$}

In this section we briefly recollect some basic results (without proofs)
concerning the Hopf algebras $\mathcal B(E)$.

Let $n \in \mathbb N^*$ and let $E \in GL(n)$.
We have already defined the algebra $\B(E)$. It was introduced by 
M. Dubois-Violette and G. Launer in \cite{[DVL]}. The 
following result is taken from \cite{[DVL]}:

\begin{prop}
The algebra $\mathcal B(E)$ admits a Hopf algebra structure, with
comultiplication $\Delta$ defined by $\Delta(a_{ij})
= \sum_{k=1}^n a_{ik} \otimes a_{kj}$, $1 \leq i,j \leq n$,
with counit $\varepsilon$ defined by $\varepsilon(a_{ij}) = \delta_{ij}$,
$1 \leq i,j \leq n$, and with antipode $S$ defined on the matrix 
$a = (a_{ij})$ by
$S(a) = E^{-1}{^t \! a} E$. $\square$
\end{prop} 

The Hopf algebra $\B(E)$ was defined in \cite{[DVL]} as the function
algebra on the quantum group of a bilinear form associated with $E$.
This is explained by the following result, which was not 
explicitly stated in \cite{[DVL]}, but was clearly implicit in that paper:

\begin{prop}
i) Consider the vector space $V = k^n$ with its canonical basis 
$(e_i)_{1 \leq i \leq n}$. Endow $V$ with the $\B(E)$-comodule structure 
defined by $\alpha(e_i) = \sum_{j=1}^n e_j \otimes a_{ji}$, 
$1 \leq i \leq n$. Then the linear map $\beta : V \otimes V \longrightarrow k$
defined by $\beta (e_i \otimes e_j) = \lambda_{ij}$, ${1 \leq i,j \leq n}$,
where $E = (\lambda_{ij})$, is a $\B(E)$-comodule morphism.

\noindent
ii) Let $A$ be a Hopf algebra and let $V$ be a finite-dimensional 
$A$-comodule of dimension $n$. 
Let $\beta : V \otimes V \longrightarrow k$ be an $A$-comodule morphism 
such that the associate bilinear form is non-degenerate.
Then there exists $E \in GL(n)$ such that $V$ is a $\B(E)$-comodule,
such that $\beta$ is a $\B(E)$-comodule morphism, and there exists 
a unique Hopf algebra morphism $\phi : \B(E) \longrightarrow A$
such that $({\rm id}_V \otimes \phi) \circ \alpha = \alpha'$, where $\alpha$
and $\alpha'$ denote the coactions on $V$ of $\B(E)$ and $A$ respectively. 
$\square$
\end{prop}

The next result was also known in \cite{[DVL]}. It will be generalized 
at the Hopf-Galois extension level in the next section.

\begin{prop}
Let $E,P \in GL(n)$. Then the Hopf algebras 
$\B(E)$ an $\B({^t \! P}EP)$ are isomorphic. $\square$
\end{prop}

We end the section by connecting the Hopf algebras $\mathcal B(E)$
with the Hopf algebra  $\qsl$. Let $q \in k^*$ and let 
$E_q = \left(\begin{array}{cc} 0 & 1 \\
                          -q^{-1} & 0\\
       \end{array} \right) \in GL(2)$.
Then it is a straightforward computation to check that
$\B(E_q) = \qsl$ (with the definition of \cite{[Ka]} for
$\qsl$).     

\section{The Hopf bigalois extensions $\mathcal B(E,F)$}

In order to prove Theorem 1.1, we introduce appropriate 
Hopf bigalois extensions. By Schauenburg's Theorem 5.5 in \cite{[Sc]},
it is  equivalent to construct Hopf bigalois extensions
and monoidal equivalences between comodule categories.
Let us first recall the language of Galois extensions for
Hopf algebras (see \cite{[Mo]} for a general perspective). 

Let $A$ be a Hopf algebra. A left $A$-Galois extension (of $k$)
is a non-zero left $A$-comodule algebra $Z$ such that the linear map
$\kappa_l$ defined by the composition
\begin{equation*}
\begin{CD}
\kappa_l : Z \otimes Z @>\alpha \otimes 1_Z>>
A \otimes Z \otimes Z @>1_A \otimes m_Z>> A \otimes Z
\end{CD}
\end{equation*}
where $\alpha$ is the coaction of $A$ and $m_Z$ is the multiplication
of $Z$, is bijective.

Similarly, a right $A$-Galois extension 
is a non-zero right $A$-comodule algebra $Z$ such that the linear map
$\kappa_r$ defined by the composition
\begin{equation*}
\begin{CD}
\kappa_r : Z \otimes Z @>1_Z \otimes \beta>>
Z \otimes Z \otimes A @>m_Z \otimes 1_A>> Z \otimes A
\end{CD}
\end{equation*}
where $\beta$ is the coaction of $A$, is bijective.

Let $A$ and $B$ be Hopf algebras. An algebra $Z$ 
is said to be an $A$-$B$-bigalois extension \cite{[Sc]} if $Z$
is both a left $A$-Galois extension and a right $B$-Galois extension,
and if $Z$ is an $A$-$B$-bicomodule.
By Theorem 5.5 in \cite{[Sc]}, there exists a monoidal
equivalence between the categories Comod$(A)$ and Comod$(B)$ if and only if
there exists an $A$-$B$-bigalois extension.

\medskip

The following definition is a natural generalization
of the definition of the algebras $\B(E)$: 

\begin{defi}
Let $E \in GL(m)$ and let $F \in GL(n)$.
The algebra $\mathcal B(E,F)$ is the universal algebra with generators
$z_{ij}$, $1\leq i \leq m,1\leq j \leq n$
and satisfying the relations
$$F^{-\!1} {^t \! z} E z = I_n \ ; \   z F^{-\!1} {^t \! z} E = I_m.$$
\end{defi}

We have $\B(E,E) = \B(E)$. 
Let us first prove a generalization of Proposition 2.3:

\begin{prop}
Let $E,P \in GL(m)$ and let $F,Q \in GL(n)$. Then the algebras 
$\B(E,F)$ and $\B({^t \! P} E P, {^t \! Q}F Q)$ are isomorphic.  
\end{prop}

\noindent
\textbf{Proof}.
Let us denote by  $y_{ij}$, $1 \leq i \leq m$, $1 \leq j \leq n$, 
the generators of $\B({^t \! P} E P, {^t \! Q}F Q)$.  
Then the relations
$$(Q^{-1} F^{-1} {^t \! Q}^{-1}) {^t \! y} ({^t \! P} E P) y = I_n
\quad {\rm and} \quad 
y (Q^{-1} F^{-1} {^t \! Q}^{-1}) {^t \! y} ({^t \! P} E P) = I_m$$
 ensure that we have an algebra morphism
$\psi : \B(E,F) \longrightarrow \B({^t \! P} E P, {^t \! Q}F Q)$ defined
by $\psi(z)= Py Q^{-1}$.
The inverse map is then defined by $\psi^{-1}(y) = P^{-1} z Q$.
$\square$

\bigskip

The algebras $\B(E,F)$ are natural candidates
to be $\B(E)$-$\B(F)$-bigalois extensions. We now define several structural
maps. Let us first fix some notations.
The generators of $\mathcal B(E)$ are denoted by 
$a_{ij}$, $1 \leq i,j \leq m$ ; the generators of $\B(F)$ are denoted 
by $b_{ij}$, $1 \leq i,j \leq n$ ; the generators of $\B(E,F)$
are denoted by $z_{ij}$, $1 \leq i \leq m$, $1 \leq j \leq n$ ;
the generators of  $\B(F,E)$ are denoted by 
$y_{ij}$, $1 \leq i \leq n$, $1 \leq j \leq m$.

The reader will easily check that the algebra morphisms described below
are well-defined, and that they are coassociative.

\medskip

\noindent
$\bullet$
The algebra morphism
$\alpha : \mathcal B(E,F) \longrightarrow 
\mathcal B(E) \otimes \mathcal B(E,F)$ defined by
$$\alpha(z_{ij}) = \sum_{k=1}^m a_{ik} \otimes z_{kj},
\ 1\leq i \leq m,1\leq j \leq n,$$
endows $\B(E,F)$ with a left $\B(E)$-comodule algebra structure.

\medskip

\noindent
$\bullet$
Similarly, the algebra morphism
$\beta : \mathcal B(E,F) \longrightarrow
\mathcal B(E,F) \otimes \mathcal B(F)$ defined by
 $$\beta(z_{ij}) = \sum_{k=1}^n z_{ik} \otimes b_{kj}, \
1\leq i \leq m,1\leq j \leq n,$$
endows $\B(E,F)$ with a right $\B(F)$-comodule algebra structure.
It is clear that $\B(E,F)$ is a $\B(E)$-$\B(F)$-bicomodule. 

\medskip

We need several other algebra morphisms to prove that the 
maps $\kappa_l$ and $\kappa_r$ are bijective. Once again it is 
straightforward to check that the algebra morphisms considered below are 
well-defined.

\medskip

\noindent
$ \bullet$ We have an algebra morphism
$\phi : \B(F,E) \longrightarrow \B(E,F)^{\rm op}$ defined 
by $$\phi(y) = F^{-1} {^t \! z} E.$$

\medskip

\noindent
$ \bullet$ We have an algebra morphism
$\gamma_1 : \B(E) \longrightarrow \B(E,F) \otimes \B(F,E)$
defined by
$$\gamma_1(a_{ij}) = \sum_{k=1}^n z_{ik} \otimes y_{kj}, \
1\leq i,j \leq m.$$
Similarly we have an algebra morphism
$\gamma_2 : B(F) \longrightarrow \B(F,E) \otimes \B(E,F)$
defined by 
$$\gamma_2(b_{ij}) = \sum_{k=1}^m y_{ik} \otimes z_{kj}, \
1\leq i,j \leq n.$$

We have introduced all the ingredients to prove the following
result:

\begin{prop}
Assume that $\mathcal B(E,F) \not = 0$. Then $\mathcal B(E,F)$ is a 
$\mathcal B(E)$-$\mathcal B(F)$-bigalois extension.
\end{prop}

\noindent
\textbf{Proof}. 
Let $\eta_l$ 
be the unique linear map such that the following diagram commutes
\begin{equation*}
\begin{CD}
\B(E) \otimes \B(E,F) 
@>\eta_l>> \B(E,F) \otimes \B(E,F) \\
@V\gamma_1 \otimes 1VV
@A1 \otimes mAA \\
\B(E,F) \otimes \B(F,E) \otimes \B(E,F) @>1 \otimes \phi \otimes 1>> 
\B(E,F) \otimes \B(E,F) \otimes \B(E,F) 
\end{CD}
\end{equation*}
and similarly, let $\eta_r$ be the unique linear map such that
the following diagram commutes
\begin{equation*}
\begin{CD}
\B(E,F) \otimes \B(F) 
@>\eta_r>> \B(E,F) \otimes \B(E,F) \\
@V1 \otimes \gamma_2VV
@Am \otimes 1AA \\
\B(E,F) \otimes \B(F,E) \otimes \B(E,F) @>1 \otimes \phi \otimes 1>> 
\B(E,F) \otimes \B(E,F) \otimes \B(E,F) 
\end{CD}
\end{equation*}
Now let us note the following identities: 
$$\sum_{k=1}^m \phi(y_{lk}) z_{kj} = \delta_{lj}, \quad
1 \leq l,j \leq n, \quad
{\rm and} \quad
\sum_{k=1}^n z_{lk} \phi(y_{kj}) = \delta_{lj}, \quad
1 \leq l,j \leq m.$$
Let $x \in \B(E,F)$, let $i,k \in \{1, \ldots, m \}$
and let $j,l \in \{1, \ldots, n \}$. Then using the previous identities,
it is immediate to check that
$$\eta_l \circ \kappa_l(z_{ij} \otimes x) =
z_{ij} \otimes x \ , \ 
\kappa_l \circ \eta_l(a_{ik} \otimes x) =  a_{ik} \otimes x,$$
and 
$$\eta_r \circ \kappa_r(x \otimes z_{ij}) =
x \otimes z_{ij} \ , \ 
\kappa_r \circ \eta_r(x  \otimes b_{jl}) =  x \otimes b_{jl}.$$
Now using the facts that $\gamma_1$, $\gamma_2$ and $\phi$ are algebras 
morphisms and that the elements considered in these equations are
generators of the corresponding algebras, 
it is not hard to see  that $\eta_l$ and $\eta_r$ are
inverse isomorphisms of $\kappa_l$ and $\kappa_r$ respectively.
$\square$

\bigskip

We have now to determine when the algebra $\B(E,F)$ is 
a non-zero algebra. This is done in the next Proposition.
Let $q \in k^*$: we use the matrix $E_q$ introduced in the
previous section.

\begin{prop}
Let $q \in k^*$ and let $F \in GL(n)$, $n \geq 2$. Assume that 
$q^2 + {\rm tr}(E ^t \!  E^{-1})q +1 = 0$. Then
$\mathcal B(E_q,F)$ is a non-zero algebra.
\end{prop}

The proof of Proposition 3.4 requires some work and 
will be done in the next section. Taking this result for
guaranteed, we can prove Theorem 1.1.
Indeed by Proposition 3.4, for $q \in k^*$ satisfying
$q^2 + {\rm tr}(E ^t \!  E^{-1})q +1 = 0$, the algebra 
$\mathcal B(E_q,F)$ is non-zero algebra. Hence 
by Proposition 3.3 $\mathcal B(E_q,F)$ is a 
$\B(E_q)$-$\B(F)$-bigalois extension. We can use Theorem 5.5 in
\cite{[Sc]} : we have an equivalence of monoidal categories
$${\rm Comod}(\mathcal B(E_q)) \cong^{\otimes} 
{\rm Comod}(\B(F))$$
and since $\B(E_q) = \qsl$, the proof of Theorem 1.1 is complete.
Let us note that this monoidal equivalence also induces 
a monoidal equivalence between the categories of  finite-dimensional
comodules. 

\medskip

\noindent
\textbf{Remark}. The algebra $\mathcal B(E,F)$ is non-zero
when ${\rm tr}(E ^t \!  E^{-1}) = {\rm tr}(F ^t \!  F^{-1})$.
We will prove this fact at the end of section 4.
 
\section{Proof of Proposition 3.4}

This section is devoted to the proof of Proposition 3.4.
Our strategy is the following one. We write a convenient presentation
for $\B(E_q,F)$ and use  Bergman's diamond lemma \cite{[Be]}
to get several linearly independent elements: this will
clearly imply that $\B(E_q,F)$ is a non-zero vector space.

Let $F= (\alpha_{ij}) \in GL(n)$ ($n\geq 2$) with inverse 
$F^{-1} = (\beta_{ij})$, and let $q \in k^*$ be such that
$q^2 + {\rm tr}(E ^t \!  E^{-1})q +1 = 0$. This equation may be rewritten in 
the most convenient form:
$${\rm tr}(E {^t \! E}^{-1}) = \sum_{i,j} \alpha_{ij} \beta_{ij}
= -q - q^{-1}.$$
We would like to be able to assume that $\beta_{nn} = 0$.
This will avoid some overlap ambiguities in the presentation of
$\mathcal B(E_q,F)$. The following elementary
lemma will be useful for this purpose:

\begin{lemm}
Let $M = (M_{ij}) \in GL(n)$ ($n \geq 2$). 
Then there exists $P \in GL(n)$ such that
$({^t \! P}MP)_{nn} = 0$.
\end{lemm}

\noindent
\textbf{Proof}. First assume that
$M_{11} = 0$. Let $P = \sum_{i=1}^n E_{n-i+1,i} \in GL(n)$, where
the $E_{ij}$'s denote the standard elementary matrices. Then
$({^t P}MP)_{nn} = M_{11} = 0$ 
Now assume that $M_{11} \not = 0$ and $M_{nn} \not = 0$.
Let $\lambda \in k^{*}$ be such that 
$\lambda^2 M_{nn} + (M_{n1} + M_{1n}) \lambda + M_{11} = 0$.
Let $P = \sum_{i=1}^{n-1} E_{ii} + \lambda E_{nn} + E_{1n} \in GL(n)$.
Then we have $({^t \! P}MP)_{nn} = 0$. $\square$

\medskip

Now take $M = F^{-1}$ and pick $P \in GL(n)$ as in Lemma 4.1:
$({^t \! P}F^{-1}P)_{nn} = 0$. Then by Proposition 3.2, the algebras
$\B(E_q,F)$ and $\B(E_q, P^{-1} F {^t \! P}^{-1})$ are isomorphic and we have
$(P^{-1} F{^t \! P}^{-1})^{-1}_{nn} = ({^t \! P}F^{-1} P)_{nn} = 0$. 
Thus we may assume without loss of generality that $\beta_{nn} = 0$.
 
\medskip

Let us now study in detail the algebra $\B(E_q,F)$: it is the universal
algebra with generators $z_{ij}$, $1\leq i \leq 2$, $1 \leq j \leq n$, and 
relations:
$${^t \! z} E_q z = F \quad {\rm and} \quad   
z F^{-\!1} {^t \! z}  = E_q^{-1}.$$
Let us write these relations explicitly:
$$z_{1i} z_{2j} - q^{-1} z_{2i} z_{1j} =  \alpha_{ij}, \quad
1 \leq i,j \leq n \ ;$$
$$\sum_{i,j=1}^n \beta_{ij} z_{1i} z_{1j} = 0
= \sum_{i,j = 1}^n \beta_{ij} z_{2i} z_{2j} \ ;$$ 
$$\sum_{i,j=1}^n \beta_{ij} z_{1i} z_{2j} = -q \quad ; \quad
\sum_{i,j=1}^n \beta_{ij} z_{2i }z_{1j} = 1.$$
Multiplying the first relation by $\beta_{ij}$, summing over
$i$ and $j$, using the third relation and the identity
${\rm tr}(E {^t \! E}^{-1}) = -q -q^{-1}$, we see that the last
relation is a consequence of the other ones.

Let us order the set $\{1, \ldots, n \} \times \{1, \ldots, n \}$
lexicographically. Take $(u,v)$ the maximal element such that
$\beta_{uv} \not = 0$. Since $\beta_{nn} = 0$, we have
$v < n$ and since the matrix $(\beta_{ij})$ is invertible, we have
$u=n$. We see now that $\B(E_q,F)$ is the universal algebra
with generators $z_{1i}$, $z_{2i}$, $1 \leq i \leq n$,
and satisfying the relations:
$$\left\lbrace
\begin{array}{l}
z_{2i} z_{1j} = q(z_{1i} z_{2j} - \alpha_{ij}) \ , \quad 
1 \leq i,j \leq n  \quad ; \\
\\
z_{1n} z_{1v} = \beta_{nv}^{-1}
(-\sum_{(i,j) < (n,v)} \beta_{ij} z_{1i} z_{1j}) \quad ; \\
\\
z_{1n} z_{2v} = \beta_{nv}^{-1}
(-q -\sum_{(i,j) < (n,v)} \beta_{ij} z_{1i} z_{2j}) \quad ; \\
\\
z_{2n} z_{2v} = \beta_{nv}^{-1}
(-\sum_{(i,j) < (n,v)} \beta_{ij} z_{2i} z_{2j}). 
\end{array}
\right.$$
We have now a nice presentation to use the diamond lemma \cite{[Be]}.
We freely use the techniques and definitions involved in the diamond lemma, 
and in particular the simplified exposition in the book 
\cite{[KS]} (although there are a few misprints there).
We endow the set $\{z_{ij}$, $1\leq i \leq 2$, $1 \leq j \leq n \}$, with
the order induced by the lexicographic order on the
set $\{1, 2 \} \times \{1, \ldots, n \}$, and we order
the set of monomials lexicographically. It is clear that the
presentation above is compatible with this order. It is also clear
that there are no inclusion ambiguities.
There are exactly the following overlap ambiguities:
$$(z_{2i}z_{1n},z_{1n}z_{1v}) \ , \quad 1\leq i \leq n \quad ; \quad
(z_{2i}z_{1n},z_{1n}z_{2v}) \ , \quad 1\leq i \leq n \quad ;$$
$$(z_{1n}z_{2v},z_{2v}z_{1j}) \ , \quad 1\leq j \leq n \quad ; \quad
(z_{2n}z_{2v},z_{2v}z_{1j}) \ , \quad 1\leq j \leq n.$$
Note that that if we had $v=n$ ($\beta_{nn} \not = 0$), there
would be more ambiguities.
We must check now that these ambiguities are resolvable.
Let us first note some preliminary identities ($i,j \in \{1, \ldots, n\}$):
$$
\sum_{(k,l)<(n,v)}   \alpha_{ik} \beta_{kl}z_{1l} =
\sum_{k,l} \alpha_{ik} \beta_{kl} z_{1l}
-\alpha_{in} \beta_{nv} z_{1v}
= z_{1i} - \alpha_{in} \beta_{nv}z_{1v}.\leqno(1)$$
Similarly, we have:
$$ 
\sum_{(k,l)<(n,v)}  \alpha_{ik} \beta_{kl}z_{2l} =
z_{2i} - \alpha_{in} \beta_{nv}z_{2v} \ ; \leqno(2)
$$
$$
\sum_{(k,l)<(n,v)}  \alpha_{lj} \beta_{kl}z_{1k}
= z_{1j} - \alpha_{vj} \beta_{nv} z_{1n} \ ; \leqno(3)
$$
$$
\sum_{(k,l)<(n,v)}  \alpha_{lj} \beta_{kl}z_{2k}
= z_{2j} - \alpha_{vj} \beta_{nv} z_{2n} \ ; \leqno(4)
$$
Let us check now that the overlap ambiguities of the first family 
are resolvable.
As usual the symbol ``$\rightarrow$'' means that we perform 
a reduction. 
Let $i \in \{1, \ldots, n\}$. We have:
\begin{align*}
& q (z_{1i}z_{2n} - \alpha_{in})z_{1v} =
q(z_{1i}z_{2n}z_{1v} - \alpha_{in}z_{1v})
\rightarrow
q \left(z_{1i}(q(z_{1n}z_{2v} - \alpha_{nv})) - \alpha_{in} z_{1v} \right) = \\
& = q(q z_{1i} z_{1n}z_{2v} - q \alpha_{nv} z_{1i} - \alpha_{in}z_{1v}) \\
& \rightarrow
q\left(q z_{1i} \left(\beta_{nv}^{-1}(-q 
-\sum_{(k,l) < (n,v)} \beta_{kl} z_{1k} z_{2l})\right)
 - q \alpha_{nv} z_{1i} - \alpha_{in} z_{1v} \right) \\
& = q\left( -q^2  \beta_{nv}^{-1} z_{1i} 
 - q \beta_{nv}^{-1}( 
\sum_{(k,l) < (n,v)} \beta_{kl} z_{1i} z_{1k} z_{2l})
 - q \alpha_{nv} z_{1i} - \alpha_{in} z_{1v}\right).
\end{align*}
On the other hand we have:
\begin{align*}
& z_{2i} 
\left(\beta_{nv}^{-1}(-\sum_{(k,l) < (n,v)} \beta_{kl} z_{1k} z_{1l})\right)
= -\beta_{nv}^{-1}\left(
 \sum_{(k,l) < (n,v)} \beta_{kl} z_{2i} z_{1k} z_{1l} \right) \\
& \rightarrow \ 
-\beta_{nv}^{-1}
\left(\sum_{(k,l) < (n,v)} \beta_{kl} 
(q(z_{1i} z_{2k}-\alpha_{ik})) z_{1l} \right) \\
& = -\beta_{nv}^{-1}
\left(\sum_{(k,l) < (n,v)} q \beta_{kl} z_{1i} z_{2k} z_{1l}
- q \sum_{(k,l) < (n,v)} \alpha_{ik} \beta_{kl} z_{1l} \right) \\
& = -\beta_{nv}^{-1} 
\left(\sum_{(k,l) < (n,v)} q \beta_{kl} z_{1i} z_{2k} z_{1l}
-q (z_{1i} - \alpha_{in} \beta_{nv} z_{1v}) \right) \ {\rm by} \ (1) \\
& \rightarrow  -\beta_{nv}^{-1} 
\left(\sum_{(k,l) < (n,v)} q \beta_{kl} z_{1i} 
(q(z_{1k}z_{2l} - \alpha_{kl}))
-q (z_{1i} - \alpha_{in} \beta_{nv} z_{1v}) \right) \\
& = -\beta_{nv}^{-1} 
\left(\sum_{(k,l) < (n,v)} 
q^2 \beta_{kl} z_{1i} z_{1k}z_{2l} 
-q^2(\sum_{(k,l) < (n,v)} \alpha_{kl} \beta_{kl} z_{1i})
-q (z_{1i} - \alpha_{in} \beta_{nv} z_{1v}) \right) \\
& = -\beta_{nv}^{-1}q 
\left(\sum_{(k,l) < (n,v)} 
q \beta_{kl} z_{1i} z_{1k}z_{2l} 
+q\alpha_{nv}\beta_{nv}z_{1i} +q^2z_{1i} +z_{1i}
- (z_{1i} - \alpha_{in} \beta_{nv} z_{1v}) \right) \\
& = q\left( -q^2  \beta_{nv}^{-1} z_{1i} 
 - q \beta_{nv}^{-1}( 
\sum_{(k,l) < (n,v)} \beta_{kl} z_{1i}z_{1k} z_{2l})
 - q \alpha_{nv} z_{1i} - \alpha_{in} z_{1v}\right).
\end{align*} 
Hence the overlap ambiguities of the first family are resolvable.
Let us now study the second family of ambiguities. We have:
\begin{align*}
& q (z_{1i}z_{2n} - \alpha_{in})z_{2v} =
q(z_{1i}z_{2n}z_{2v} - \alpha_{in}z_{2v})
\rightarrow
q \left(-z_{1i} \beta_{nv}^{-1} (\sum_{(k,l)< (n,v)} \beta_{kl}z_{2k}z_{2l})
 - \alpha_{in} z_{2v} \right) \\ 
& = q \left(-\beta_{nv}^{-1} (\sum_{(k,l)< (n,v)}
\beta_{kl}z_{1i} z_{2k} z_{2l})
 - \alpha_{in} z_{2v} \right).
\end{align*}
On the other hand we have:
\begin{align*}
& z_{2i}\left(\beta_{nv}^{-1}
(-q-\sum_{(k,l)< (n,v)} \beta_{kl}z_{1k}z_{2l})\right)
= \beta_{nv}^{-1}
\left(-\sum_{(k,l)< (n,v)} \beta_{kl} z_{2i}z_{1k}z_{2l}
-qz_{2i}\right) \\
& \rightarrow
\beta_{nv}^{-1}
\left(-\sum_{(k,l)< (n,v)} \beta_{kl} q (z_{1i}z_{2k} - \alpha_{ik})z_{2l}
-qz_{2i}\right) \\
& = \beta_{nv}^{-1}q
\left(- \sum_{(k,l)< (n,v)} \beta_{kl} z_{1i}z_{2k}z_{2l}
+ \sum_{(k,l)< (n,v)} \alpha_{ik} \beta_{kl} z_{2l}
-z_{2i}\right) \\
& = \beta_{nv}^{-1}q
\left(-\sum_{(k,l)< (n,v)} \beta_{kl} z_{1i}z_{2k}z_{2l}
+ z_{2i} - \alpha_{in} \beta_{nv}z_{2v}
-z_{2i}\right) \quad {\rm by} \ (2) \\
& = q
\left(-\beta_{nv}^{-1}
(\sum_{(k,l)< (n,v)} \beta_{kl} z_{1i}z_{2k}z_{2l})
- \alpha_{in} z_{2v}\right).
\end{align*}
Hence the ambiguities in the second family are resolvable.
The resolvability of the ambiguities of the third and fourth
families is shown is the same way, using the identities 
(3) and (4) respectively. This is left to the reader.
Since all ambiguities are resolvable and our order is compatible with
the presentation, we can use the diamond lemma \cite{[Be]}:
the set of reduced monomials (i.e. those invariant under all
reductions) is a basis of $\B(E_q,F)$. In particular 
the reduced monomials $z_{1i}$, $1 \leq i \leq n$, are linearly independent.
This shows that $\B(E_q,F)$ is a non-zero vector space, and  concludes the
proof of Proposition 3.4. $\square$

\bigskip

Let $E \in GL(m)$ and $F \in GL(n)$.
Let us prove now that $\mathcal B(E,F)$ is non-zero
when ${\rm tr}(E ^t \!  E^{-1}) = {\rm tr}(F ^t \!  F^{-1})$.

Let $q \in k^*$ be such that
$q^2 + {\rm tr}(E ^t \!  E^{-1})q +1 = 0$.
Similarly as in Section 3, we have an algebra morphism 
$\delta : \mathcal B(E,F) \longrightarrow 
\mathcal B(E, E_q) \otimes \mathcal B(E_q,F)$ defined by
$$\delta(z_{ij}) = \sum_{k=1}^2 v_{ik} \otimes w_{kj}, \
1\leq i \leq m,1\leq j \leq n,$$
where $(v_{ik})$ and $(w_{kj})$ denote the generators of 
$\B(E,E_q)$ and $\B(E_q,F)$ respectively. Then by the proof of Proposition
3.4, the elements $w_{kj}$ are linearly independent, and
since the algebra morphism $\phi : \B(E,E_q) \longrightarrow \B(E_q,E)^{\rm op}$
of Section 3 is an isomorphism, the elements $v_{ik}$ are also linearly 
independent. Hence $\delta(z_{ij}) \not = 0$, and it follows that
$\B(E,F)$ is  a non-zero algebra.

\section{$SL(2)$-deformations}

In this section $k$ will be an algebraically closed field of characteristic
zero. This section is essentially devoted to the proof of
Theorem 1.2. We also determine the isomorphic classification of the
Hopf algebras $\B(E)$.

Let us first recall the concept of representation semi-ring of 
a Hopf algebra. Let $A$ be a cosemisimple Hopf algebra. The representation ring
of $A$ is defined to be the Grothendieck group of the category
Comod$_f(A)$ : $R(A) = K_0({\rm Comod}_f(A))$. 
It is a free abelian
group with a basis formed by the isomorphism classes of
simple (irreducible) comodules.
The monoidal
structure of Comod$_f(A)$ induces a ring structure on $R(A)$.
The isomorphism class of a finite-dimensional $A$-comodule
$V$ is denoted by $[V]$. The representation semi-ring of 
$A$ is now defined to be
$$R^+(A) = \{ \sum_i a_i [V_i] \in R(A), \ a_i \geq 0, 
V_i \in {\rm Comod}_f(A) \}.$$
Let $B$ be another cosemisimple Hopf algebra and 
let $f : A \longrightarrow B$ be a Hopf algebra morphism. Then $f$ induces 
a monoidal 
functor $f_* : {\rm Comod}_f(A) \rightarrow {\rm Comod}_f(B)$, and hence
a semi-ring morphism $f_* :R^+(A) \longrightarrow R^+(B)$. 
It is not difficult to see that a semi-ring isomorphism
$R^+(A) \cong R^+(B)$ induces a bijective correspondence
(that preserves tensor products) between the isomorphism classes
of simple comodules of $A$ and $B$.

\smallskip

Let $G$ be a reductive algebraic group. It is classical to say
that a cosemisimple Hopf algebra $A$ is a $G$-deformation if 
one has a semi-ring isomorphism $R^+(A) \cong R^+(\mathcal O(G))$.
Hence Theorem 2.1 classifies $SL(2)$-deformations.

\medskip

Let us state a useful folk-known result. We include a sketch of
proof for the sake of completeness.

\begin{lemm}
Let $A$ and $B$ be cosemisimple 
Hopf algebras and let $f : A \longrightarrow B$ be 
a Hopf algebra morphism inducing a semi-ring isomorphism
$R^+(A) \cong R^+(B)$. Then $f$ is a Hopf algebra isomorphism.
\end{lemm}

\noindent
\textbf{Proof}. let $f_*  : {\rm Comod}_f(A) \rightarrow {\rm Comod}_f(B)$
be the induced functor. Since $f$ induces an isomorphism
$R^+(A) \cong R^+(B)$, the functor $f_*$ transforms simple objects
of ${\rm Comod}_f(A)$ into simple objects of  ${\rm Comod}_f(A)$ and 
hence is an equivalence of categories (the categories 
$ {\rm Comod}_f(A)$ and $ {\rm Comod}_f(B)$ are semisimple).
Then $f$ is an isomorphism by Tannaka-Krein type reconstruction theorems
(see e.g. \cite{[JS]}). $\square$.

\bigskip

We now recall the representation theory of
$\qsl$. Our reference, especially for the root of unity case,
 will be the paper \cite{[KP]} of 
P. Kondratowicz and P. Podle\'s.

Let $q \in k^*$. We say that $q$ is generic if $q$ is not a root of unity
or if $q \in \{\pm 1 \}$.

\noindent
$\bullet$ Let us first assume that
$q$ is generic. Then $\qsl$ is cosemisimple and has a 
family of non-isomorphic simple comodules $(U_n)_{n \in \N}$ such that
$$U_0 = k, \quad U_n \otimes U_1 \cong U_1 \otimes U_n
\cong U_{n-1} \oplus U_{n+1},
\quad \dim(U_n) = n+1, \ \forall n \in \N^*.$$
Furthermore, any simple $\qsl$-comodule is isomorphic to one of the
comodules $U_n$.

\medskip

\noindent
$\bullet$ Now assume that $q$ is not generic. Let $N \geq 3$ be the order
of $q$. Then we let
\begin{equation*}
N_0 = 
\begin{cases} N & \text{if $N$ is odd},\\
N_0/2 & \text{if $N$ is even}.
\end{cases}
\end{equation*}
Then $\qsl$ is not cosemisimple. There exists families
$\{V_n, \ n \in \N \}$, $\{U_n, \ n= 0, \ldots, N_0-1 \}$ of non-isomorphic
simple comodules (except for $n = 0$ where $U_0 \cong k \cong V_0$),
such that
$$V_n \otimes V_1 \cong V_1 \otimes V_n
\cong V_{n-1} \oplus V_{n+1},
\quad \dim(V_n) = n+1, \ \forall n \in \N^*.$$
$$U_n \otimes U_1 \cong U_1 \otimes U_n
\cong U_{n-1} \oplus U_{n+1},
\quad \dim(U_n) = n+1, \ n = 1, \ldots, N_0-2 .$$
The comodule $U_{N_0-1} \otimes U_1$ is not semisimple.
It has a simple filtration:
$$(0)  \subset U_{N_0-2} \subset Y \subset U_{N_0-1} \otimes U_1,$$
with $U_{N_0-1} \otimes U_1 / Y \cong U_{N_0-2}$ and 
$Y/ U_{N_0-2} \cong V_1$. 

The comodules $V_n \otimes U_m  \cong U_m \otimes V_n$, 
$n \in \N$, $m = 0, \ldots , N_0-1$,
are simple, and any simple $\qsl$-comodule is isomorphic with one of
these comodules.

Finally there is another useful fact: the Hopf subalgebra of $\qsl$
generated by the matrix coefficients of the comodule $V_1$ is cosemisimple
and is isomorphic with $\mathcal O_{\pm 1}(SL(2))$. 

Let $E \in GL(m)$, $m \geq 2$. Let $n \in \N$. We denote by $U_n^E$ and $V_n^E$
the simple $\mathcal B(E)$-comodules corresponding to 
the simple $\qsl$-comodules $U_n$ and $V_n$
(for $q$ as in Theorem 1.1).

\medskip

Here is a useful lemma:

\begin{lemm}
Let $E \in GL(m)$ and $F \in GL(n)$. Let
$F : {\rm Comod}_f(\B(E)) \rightarrow {\rm Comod}_f(\B(F))$ be
an equivalence of monoidal categories. Then
$F(U_1^E) \cong U_1^F$. If $\B(E)$ is cosemisimple, then
$F(U_n^E) \cong U_n^F$, $\forall n \in \N$.
If $\B(E)$ is not cosemisimple, then
$F(U_n^E) \cong U_n^F$, $\forall n < N_0-1$, and
$F(V_n^E) \cong V_n^F$, $\forall n \in \N$.
\end{lemm}

\noindent
\textbf{Proof}. Let us first assume that $\B(E)$ is cosemisimple.
One can show by induction that  
$$U^E_k \otimes U^E_l \cong U^E_{|k-l|} \oplus U^E_{|k-l| +2}
\oplus \ldots \oplus U^E_{k+l} \ 
{\rm for} \  k,l \in \N.$$
Hence $U^E_1$ is the only simple $\B(E)$-comodule $W$ such that
$W \otimes W$ is the direct sum of two simple comodules. 
It follows that $F(U^E_1) \cong U^F_1$, and then an easy induction,
using the fusion rule 
$U^E_1 \otimes U^E_n \cong U^E_{n-1} \oplus U^E_{n+1}$,
shows that $F(U^E_n) \cong U^F_n$, $\forall n \in  \N$.

Now assume that $\B(E)$ is not cosemisimple. The 
$\B(F)$-comodule $F(V^E_1)^{\otimes k}$ must be semisimple for all $k \in \N$,
and hence we have $F(V^E_1) \cong V^F_i$ for some $i$. By the 
cosemisimple case we have $F(V^E_1) \cong V^F_1$ and 
$F(V^E_n) \cong V^F_n$, $\forall n \in \N$.
Now pick $Z$ a simple $\B(E)$-comodule such that $F(Z) \cong U^F_1$.
We have $Z \cong U^E_i \otimes V^E_j$ for some positive integers
$i \leq N_0-1$ and $j$. Hence $F(U^E_i) \otimes V_j^F \cong U^F_1$.
It is then clear that $j = 0$ and $F(U^E_i) \cong U^F_1$. Now note that since 
$F(U^E_k) \otimes  F(U^E_1) \cong F(U^E_{k-1}) \oplus F(U^E_{k+1})$, we have
$$\dim(F(U^E_1)) <  \dim(F(U^E_2)) < \ldots  < \dim(F(U^E_{N_0-1})).$$
Then if $i>1$, we have $\dim(F(U^E_1)) < \dim(U^F_1)$. But another
glance at the fusion rules shows that $U_1^F$ is the simple
$\B(F)$-comodule of the smallest dimension. Hence $i=1$ and 
$F(U^E_1) \cong U^F_1$. Another easy induction now shows that
$F(U^E_n) \cong U^F_n$, for $n \in  \{0, \ldots, N_0-1 \}$. $\square$

\bigskip

We now present the isomorphic classification of the Hopf
algebras $\B(E)$.

\begin{theo}
Let $E \in GL(m)$ and $F \in GL(n)$. The Hopf algebras
$\B(E)$ and $\B(F)$ are isomorphic if and only if
$m=n$ and there exists $M \in GL(m)$
such that $F = {\!^tM}EM$.
\end{theo}

\noindent
\textbf{Proof}. We denote by $(a_{ij})$ and $(b_{ij})$ the
respective generators of $\B(E)$ and $\B(F)$, and by $a$ and $b$ the
corresponding matrices. By the construction of the categorical
equivalence of Theorem 1.1 (see \cite{[U]} and \cite{[Sc]}),
the elements $a_{ij}$ and $b_{ij}$ are the matrix coefficients of 
the comodules $U_1^E$ and $U_1^F$ respectively.

Let $f : \B(E) \longrightarrow \B(F)$ be a Hopf algebra isomorphism
and let $f_* : {\rm Comod}_f(\B(E)) \longrightarrow 
{\rm Comod}_f(\B(F))$ be the induced equivalence of monoidal categories.
By Lemma 5.2, 
we have $f_*(U_1^E) \cong  U_1^F$. Hence $m=n$ and
there exists $P \in GL(m)$ such that $f(a) = PbP^{-1}$. 
But we must have $f(E^{-\!1} {^t \! a} E a) = I$, and hence
$S(b) = b^{-1}= ({^t \! P} E P)^{-1} {^t \! b} ({^t \! P} E P) =
F^{-1}{^t \! b}F$. Since the elements
$b_{ij}$ are linearly independent (the comodule $U_1^F$ is simple),
it follows that there exists $\lambda \in k^*$ such that
$F = \lambda {^t \! P} E P$, and we can take $M = \sqrt{\lambda} P$.

The converse assertion is Proposition 2.3. $\square$

\begin{rem}
The proof of Theorem 5.3 also shows that the automorphism group of 
the Hopf algebra $\B(E)$ is isomorphic
with the group
$G_E = \{ P \in GL(m) \ | \ {^t \!P}EP = E \}/ \{ \pm I \}$. $\square$   
\end{rem}

\noindent
{\large \textbf{Proof of Theorem 1.2.}}

\medskip

The proof follows closely some part of the proof of theorem 3.2
in \cite{[PM]}. Let $A$ be a cosemisimple Hopf algebra with
$R^+(A) \cong R^{+}(\mathcal O(SL(2))$. Let us denote by $U_n^A$, $n \in \N$, 
the simple $A$-comodules (with the same conventions as before).
We have $U_1^A \otimes U_1^A \cong k \oplus U_2^A$. Hence the
$A$-comodule $U_1^A$ is self-dual : there exists
$A$-comodule morphisms $e :  U_1^A \otimes U_1^A \longrightarrow
k$ and $\delta : k \longrightarrow U_1^A \otimes U_1^A$ such that
$(e \otimes {\rm id}) \circ ({\rm id} \otimes \delta) = {\rm id}$
and $({\rm id} \otimes e) \circ (\delta \otimes {\rm id}) = {\rm id}$.
These equations show that the bilinear form induced by $e$ is non-degenerate.
Thus by Proposition 2.2 there exists $E \in GL(m)$ (with 
$m = \dim(U_1^A)$) and a Hopf algebra morphism $f : \B(E) \longrightarrow
A$ such that  $f_*(U_1^E) = U_1^A$.

First assume that $\B(E)$ is cosemisimple. Then using the fusion rules
$U^E_1 \otimes U^E_n \cong U^E_{n-1} \oplus U^E_{n+1}$, an easy induction
shows that $f_*(U_n^E) \cong U_n^A$, $\forall n \in \N$.
Hence $f$ induces a semi-ring isomorphism 
$R^+(\B(E)) \cong R^+(A)$, and is an isomorphism by Lemma 5.1.

Now assume that $\B(E)$ is not cosemisimple. 
An induction also shows that $f_*(U_n^E) \cong U_n^A$, $\forall n \in 
\{0, \ldots N_0 -1 \}$.
Then we have
$$f_*(U_{N_0-1}^E \otimes U_1^E) \cong U_{N_0-1}^A \otimes U_1^A
\cong U_{N_0-2}^A \oplus U_{N_0}^A.$$
On the other hand using the simple filtration of 
the $\B(E)$-comodule $U_{N_0-1}^E \otimes U_1^E$, we have 
$$f_*(U_{N_0-1}^E \otimes U_1^E) \cong U_{N_0-2}^A \oplus 
f_*(V_1^E) \oplus U_{N_0-2}^A.$$
This contradicts the unicity of the decomposition a semisimple
comodule into a direct sum of simple comodules. 

Thus $\B(E)$ is cosemisimple,
any element $q \in k^*$ such that 
$q^2 + {\rm tr}(E ^t \!  E^{-1})q +1 = 0$ is generic,
and $f$ is an isomorphism. 
The last assertion in Theorem 
1.2 is a direct consequence of Theorem 5.3. $\square$

\section{CQG algebra structure on $\mathcal B(E)$}

In this section $k = \C$. We determine the possible
Hopf $*$-algebra structures  and CQG algebra structures on $\B(E)$.

Let us recall that a Hopf $*$-algebra is a Hopf algebra $A$, which
is also a $*$-algebra and such that the comultiplication is a 
$*$-homomorphism. If $a = (a_{ij}) \in M_n(A)$ is a matrix with coefficients 
in $A$, the matrix $(a_{ij}^*)$ is denoted by $\overline{a}$, while 
${^t \! \overline{a}}$, the 
transpose matrix of $\overline{a}$, is denoted by $a^*$. The matrix $a$
is said to be unitary if $a^*a = I = aa^*$.
Recall \cite{[KS]} that a Hopf $*$-algebra $A$ is said to be a CQG
algebra if for every finite-dimensional $A$-comodule with associate
matrix of coefficients $a \in M_n(A)$, there exist $K \in GL(n)$
such that the matrix $KaK^{-1}$ is unitary.
A CQG algebra may be seen as the algebra of representative functions 
on a compact quantum group.  

\begin{prop}
Let $E \in GL(m)$, and denote by $a_{ij}$, $1 \leq i,j \leq m$, the
generators of $\B(E)$.

\noindent
1) The Hopf algebra $\B(E)$ admits a Hopf $*$-algebra structure
if and only if there exists $M \in GL(m)$ such that
$$(\star) \quad {^t \! M} E^* M = E, \quad \overline{M} M =
\lambda I, \quad {\rm for} \ {\rm some} \ \lambda \in \R^*.$$
The $*$-structure of $\B(E)$ is then defined by 
$\overline{a} = M a M^{-1}$.
The corresponding Hopf $*$-algebra is denoted
by $\B(E)_M$. If $N$ is another matrix satisfying the conditions $(\star)$,
the Hopf $*$-algebras $\B(E)_M$ and $\B(E)_N$ are isomorphic if and only if
there exists $P \in GL(m)$ such that 
$$E = {^t \! P}EP, \quad  MP = \gamma \overline{P} N 
\quad {\rm for} \ {\rm some} \ \gamma \in \C^*.$$
2) Let $M \in GL(m)$ satisfying the conditions $(\star)$.
Then the Hopf $*$-algebra $\B(E)_M$ is a CQG algebra if and only if
there exist $\mu \in \C^*$ such that the matrix $\mu{^t \! M}^{-1}E$ is
positive. 
\end{prop}

\noindent
\textbf{Proof}.
We use the notations of Section 5: let $U_1^E$ be the fundamental
comodule of $\B(E)$, with  $(a_{ij})_{1 \leq i,j \leq m}$ 
as matrix coefficients.

\noindent
1) Let us first assume that $\B(E)$ admits a Hopf $*$-algebra structure. 
Then by the arguments of Lemma 5.2, we have $\overline{U_1^E} \cong U_1^E$,
where $\overline{U_1^E}$ denotes the conjugate comodule of
$\B(E)$. Hence there exists $M \in GL(m)$ such that 
$\overline{a} = M a M^{-1}$. Now we have
$\overline{E^{-\!1} {^t \! a} E a} = \overline{I} = I$, but
$$\overline{E^{-\!1} {^t \! a} E a} = 
{^t \! (}{^t \! (} \overline{Ea}) ^t \! (\overline{E^{-1} {^t \! a}})) = 
{^t \! (}{^t \! (} \overline{E}M a M^{-1}) 
^t \! (\overline{E}^{-1} {^t \! M}^{-1} {^t \! a} {^t \! M})) = 
{^t \! (}{^t \! M^{-1}} {^t \! a} {^t \! M} E^* 
M a M^{-1} E^{*^{-1}}),
$$
and hence 
$(M^{-1}  E^{*^{-1}}
{^t \! M^{-1}}) {^t \! a} ({^t \! M} E^* M) a  = I$.
By the linear independence of the elements $a_{ij}$, we have
$E = \alpha {^t \! M} E^* M$ for some $\alpha \in \C^*$,
and up to a normalization by $\sqrt{\alpha}$, we can assume that 
$E = {^t \! M} E^* M$. Now we have
$a = \overline{\overline{a}} =  \overline{Ma M^{-1}} =
\overline{M} M a M^{-1} \overline{M}^{-1}$. Hence by the linear independence
of the elements $a_{ij}$, we have $\overline{M}M = \lambda I$ for 
some $\lambda \in \R^{*}$.

Conversely, if $M \in GL(m)$ satisfies the conditions $(\star)$, it is
easy to check, using the computations already done, that one defines a 
Hopf $*$-algebra structure on $\B(E)$ by letting $\overline{a} =
M a M^{-1}$.

Let $M,N \in GL(m)$ satisfying the conditions $(\star)$, and let 
$\phi : \B(E)_M \longrightarrow \B(E)_N$ be a Hopf $*$-algebra isomorphism. 
By Theorem 5.3 and its proof there exists $P$ in $GL(m)$ such that 
$E = {^t \!P} E P$ and $\phi(a) = PaP^{-1}$. But we have
$$\phi(\overline{a}) = \phi(Ma M^{-1}) = MPaP^{-1}M^{-1}
\quad {\rm and}\quad \overline{\phi(a)}
= \overline{P} N a N^{-1} \overline{P}^{-1},$$
which implies that $MP = \gamma \overline{P} N$, for some $\gamma \in 
\C^*$. It is not difficult to prove the converse assertion using 
the above considerations.

\smallskip

\noindent
2) Assume that $\B(E)_M$ is a CQG algebra. Then there exists $K \in GL(m)$
such that the matrix $KaK^{-1}$ is unitary, i.e.
$$(KaK^{-1})^* (KaK^{-1}) = I = (KaK^{-1}) (KaK^{-1})^*,$$
and some easy computations show that there exists $\mu \in \C^*$ such
that ${^t \! M} K^*K = \mu E$, which means that $\mu {^t \! M}^{-1} E$
is a positive matrix.

Conversely, if there exists $\mu \in \C^*$ such that
$\mu {^t \! M}^{-1} E$ is a positive matrix, then there exists
$K \in GL(m)$ such that $\mu {^t \! M}^{-1}E = K^*K$. It is a direct 
computation to show that the matrix $KaK^{-1}$ is unitary. The algebra
$\B(E)$ is generated by the elements $a_{ij}$, and hence it follows from
\cite{[KS]}, Proposition 28, p. 417, that $\B(E)$ is a CQG algebra.
$\square$


\begin{thebibliography}{25}

\small{
\bibitem{[Ba1]} T. \textsc{Banica}, Th\'eorie des repr\'esentations du
groupe quantique compact libre $O(n)$, C.R. Acad. Sci. Paris, 322, S\'erie I,
241--244, 1996. 

\bibitem{[Ba2]} T. \textsc{Banica}, A reconstruction result for the
$R$-matrix quantizations of $SU(N)$, preprint 1998.

\bibitem{[Be]} G.M. \textsc{Bergman}, The diamond lemma for ring theory,
Adv. Math. 29, 178-218, 1978.

\bibitem{[Br]} A. \textsc{Brugui\`eres}, Dualit\'e tannakienne
pour les quasi-groupo\"{\i}des quantiques,
Comm. Algebra 25(3), 737-767, 1997.

\bibitem{[DVL]} M. \textsc{Dubois-Violette}, G. \textsc{Launer},
The quantum group of a non-degenerate bilinear form, Phys. Lett. B,
245(2), 175-177, 1990.

\bibitem{[Ph]} \textsc{Phung Ho Hai}, On matrix quantum groups
of type $A_n$, International  J. Math. 11(9), 1115-1146, 2000. 

\bibitem{[JS]} A. \textsc{Joyal}, R. \textsc{Street}, An introduction to
Tannaka duality and quantum groups, Lecture Notes in Math. 1488, 
Springer-Verlag, 1991, 413-492.

\bibitem{[Ka]} C. \textsc{Kassel}, \textsl{Quantum groups}, 
GTM 155, Springer, 1995.

\bibitem{[KS]} A. \textsc{Klimyk}, K. \textsc{Schm\"{u}dgen},
{\sl Quantum groups and their representations}, 
Texts and Monographs in Physics,
Springer, 1997. 

\bibitem{[KP]} P. \textsc{Kondratowicz}, P. \textsc{Podle\'s},
On representation theory of quantum $SL_q(2)$ at roots of unity,
Banach Center Publ. 40, 223-248, 1997.

\bibitem{[Mo]} S. \textsc{Montgomery}, 
\textsl{Hopf algebras and their actions on rings}, AMS, Providence, 1993.

\bibitem{[O]} C. \textsc{Ohn}, Quantum $SL(3, \mathbb C)$ with classical
representation theory, J. Algebra 213(2), 721-756, 1999.

\bibitem{[PM]} P. \textsc{Podle\'s}, E. \textsc{M\"{u}ller},
Introduction to quantum groups,
Rev. Math. Phys. 10, 511-551, 1998.

\bibitem{[Sc]} P. \textsc{Schauenburg}, Hopf bigalois extensions, Comm.
Algebra 24 (12), 3797--3825, 1996.

\bibitem{[U]} K.H. \textsc{Ulbrich}, Fiber functors on finite
dimensional comodules, Manuscripta Math. 65, 39--46, 1989.

\bibitem{[Wo]} S.L. \textsc{Woronowicz}, New quantum deformation
of $SL(2,\mathbb C)$. Hopf algebra level, Rep. Math. Phys. 30(2),
259-269.

}

\end{thebibliography}
\end{document}